\documentclass[hidelinks,onefignum,onetabnum]{siamart220329}

\usepackage{amsfonts,upquote}
\usepackage{graphicx,verbatim}
\headers{The first five years of the AAA algorithm}{Nakatsukasa, S\`ete, and Trefethen}

\title{The first five years of the AAA algorithm}

\author{Yuji Nakatsukasa\thanks{Mathematical Institute, U. of Oxford, Oxford OX2 6GG, UK}
  (\email{nakatsukasa@maths.ox.ac.uk}).
\and Olivier S\`ete\thanks{Inst.\ Math.\ u. Informatik, University of Greifswald, 17489 Greifswald, Germany}
  (\email{olivier.sete@ uni-greifswald.de}).
\and Lloyd N. Trefethen\thanks{School of Engineering and Applied Sciences, Harvard U., Cambridge, MA 02138, USA}
  (\email{trefethen@seas.harvard.edu}).}

\def\Re{\hbox{\rm Re\kern .5pt}}
\def\Im{\hbox{\rm Im\kern .5pt}}

\begin{document}

\maketitle

\begin{abstract}
The AAA algorithm, introduced in 2018,
computes best or near-best rational approximations to functions
or data on subsets of the real line or the complex plane.  It is much faster
and more robust than previous algorithms for such problems and has been used
in many applications since its appearance, including the numerical solution of
Laplace, Poisson, and biharmonic PDE problems in irregular domains.
AAA has also been extended in new directions and seems likely to be a
tool of lasting importance in the future.
\end{abstract}



The AAA algorithm is a numerical method for rational approximation
of a function $f$ on a real or complex domain.  The three of us
introduced the method in 2018~\cite{aaa} together with a computer
program {\tt aaa} in the Chebfun package~\cite{chebfun}, which also runs standalone
in MATLAB or Octave.\ \ AAA, pronounced ``triple A'' and derived from
``adaptive Antoulas-Anderson,''
is much faster and more robust than previous methods in this area,
and it is changing perceptions of how rational approximations can
be used in applications.

For a simple example, the MATLAB code segment

{\small
\begin{verbatim}


         Z = exp(2i*pi*(1:100)/100);
         r = aaa(exp(Z),Z);

\end{verbatim}
\par}

\noindent
computes an approximation to $\exp(z)$ in $100$ roots of unity with
the default relative accuracy of $10^{-13}$.  The computation takes
$0.002$~s on our laptop and delivers a degree~7 rational function
whose maximum deviation from $e^z$ on the unit disk is $2.81\times
10^{-15}$.  For another example, the code segment

{\small
\begin{verbatim}


         Z = linspace(4-50i,4+50i);
         zeta = @(z) sum((1e4:-1:1).^(-z),2);
         r = aaa(zeta(Z),Z);

\end{verbatim}
\par}

\noindent
evaluates the Riemann zeta function at 100 complex points with
$\Re(z) = 4$ and then fits this data in $0.1$~s by a rational
function of degree 37.  The phase portrait of $r$~\cite{wegert},
shown in Figure~1, closely captures that of the zeta function itself
in the striped region.  The first two zeros of $r$ in the upper
half-plane are at $0.4999999987 +14.1347251412i$ and $0.4999999987
+21.0220396409i$, matching the corresponding zeros of $\zeta(z)$
to 8 digits.  The pole of $r$ near $z=1$ has similar accuracy,
falling at $z\approx 1.0000000069 + 0.0000000009i$.

\begin{figure}
\begin{center}
\includegraphics[trim=30 140 30 10, clip, scale=.93]{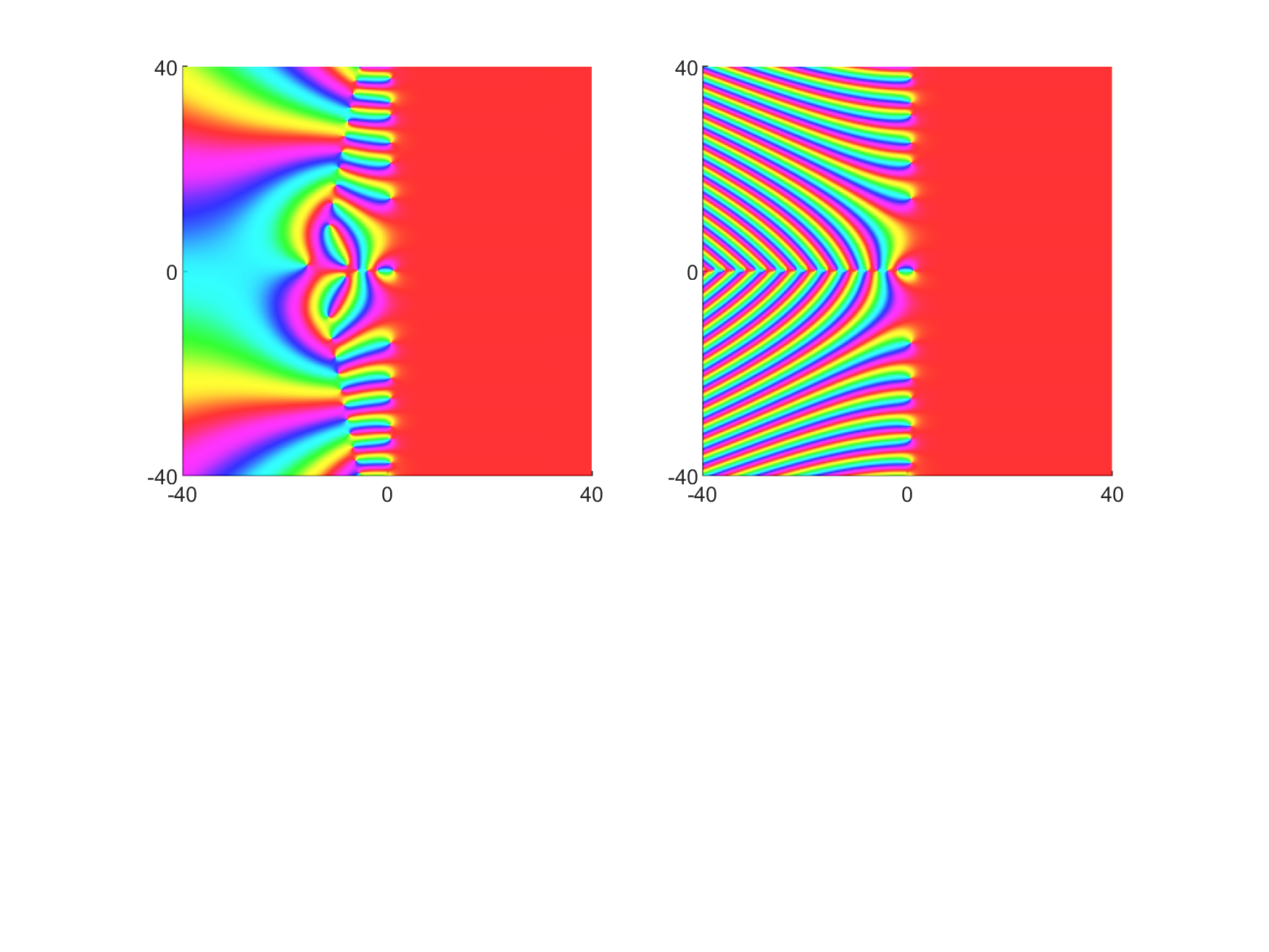}
\end{center}
\vskip -5pt
\caption{Phase portrait of the Riemann zeta function (right) compared
with that of a rational approximation based
on $100$ sample points in the right half-plane (left).
The AAA algorithm computes this approximation in $0.1$~s on a laptop.}
\end{figure}

AAA has become a standard method for many computations, with
a rapidly growing literature.  Applications to date include
analytic continuation~\cite{lustri,analcont}, interpolation of
equispaced data~\cite{equi}, Laplace problems with applications
to magnetics~\cite{costa2,costa}, conformal mapping~\cite{cm1},
Stokes flow~\cite{xue}, nonlinear dynamics~\cite{lustri,analcont}, chemical physics~\cite{kkh},
microwave tubes~\cite{tubes},
simulation of turbulence~\cite{kkw},
nonlinear eigenvalue problems~\cite{gnt,lmpv,sem},
finite element linearizations~\cite{djm}, design
of preconditioners~\cite{budisa,dpp}, model order
reduction~\cite{abg,gg,johns,pradovera}, and signal
processing~\cite{dpp,hochman,RF,vre,wdt}.  For signal processing,
AAA has been adopted as the basis of the {\tt rational} code in the MathWorks RF
Toolbox~\cite{RF}, which includes practical enhancements that our
own software lacks like the option to impose real symmetry.
There have also been generalizations to
multivariate approximation~\cite{cbg,gg,hochman,lmpv}.

The starting point of AAA is the representation of a rational function $r(z)$
not in the form $p(z)/q(z)$ with polynomials $p$ and $q$,
but by a {\em barycentric quotient}~\cite{berrut}
$$
r(z) = {n(z)\over d(z)} = \sum_{j=1}^m {w_j f(s_j)\over z-s_j} \left/ \sum_{j=1}^m {w_j\over z-s_j} \right. ,
\label{theapprox}
\eqno (1)
$$
where $s_1,\dots ,s_m$ are a set of $m$ {\em support points} and
$w_1,\dots , w_m$ are a set of nonzero {\em barycentric weights}.
This perhaps counterintuitive formula can represent a degree $m-1$
rational function in a numerically stable manner even when the zeros
and poles of $r$ are clustered near singularities, which is one of
the sources of the power of rational approximations.   The $p/q$
representation fails in such cases because $p$ and $q$ vary widely
in magnitude over the approximation domain, causing loss of accuracy
in floating point arithmetic even if they are expanded in the best
possible basis.

The AAA algorithm starts from a real or complex discrete sample set
$Z$ and a set $F$ of corresponding function values.  It consists
of an alternation between a nonlinear greedy step, in which the
next support point $z_m$ is chosen as the sample point where the
error is largest, and a linear algebra step, in which a singular
value decomposition (SVD) of a rectangular Loewner matrix $A$
is computed to determine coefficients $\{w_j\}$ to minimize the
discretized least-squares error $f(z)d(z) - n(z)$ for $z\in Z$.
The rows of $A$ correspond to sample points and the columns to
support points.  For details, including a 40-line prototype computer
code, see~\cite{aaa}.

\smallskip
{\em AAA-Lawson algorithm for minimax approximation.}  In its
original mode of operation, AAA computes a near-best as opposed to a
best approximation (by which we mean a {\em minimax} approximation,
minimizing the $L^\infty$ error, as is common in approximation theory).
However, with a second phase of iteratively reweighted least-squares
(IRLS) iteration, it is possible in most cases to improve
near-best to best, as shown in 2020~\cite{aaalawson}.  We call this
the {\em AAA-Lawson method,} since the iteration is a nonlinear
barycentric analogue of the IRLS method introduced by Lawson in
1961 for linear approximation problems~\cite{lawson}.  This is the
only known robust fast algorithm for rational best approximation
where the function or the domain or both are complex.  (For real
approximation on a real interval, there is the Remez algorithm,
which for numerical stability must also be based on the barycentric
representation~\cite{filip}.)  An illustration is given in Figure~2.
AAA-Lawson is a valuable tool for approximation theory explorations,
but we do not recommend minimax approximation as the starting point
for most applications, since standard AAA is 
more robust, especially for accuracies within a few
orders of magnitude of machine precision.

\begin{figure}
\begin{center}
\includegraphics[trim=20 132 30 0, clip, scale=.85]{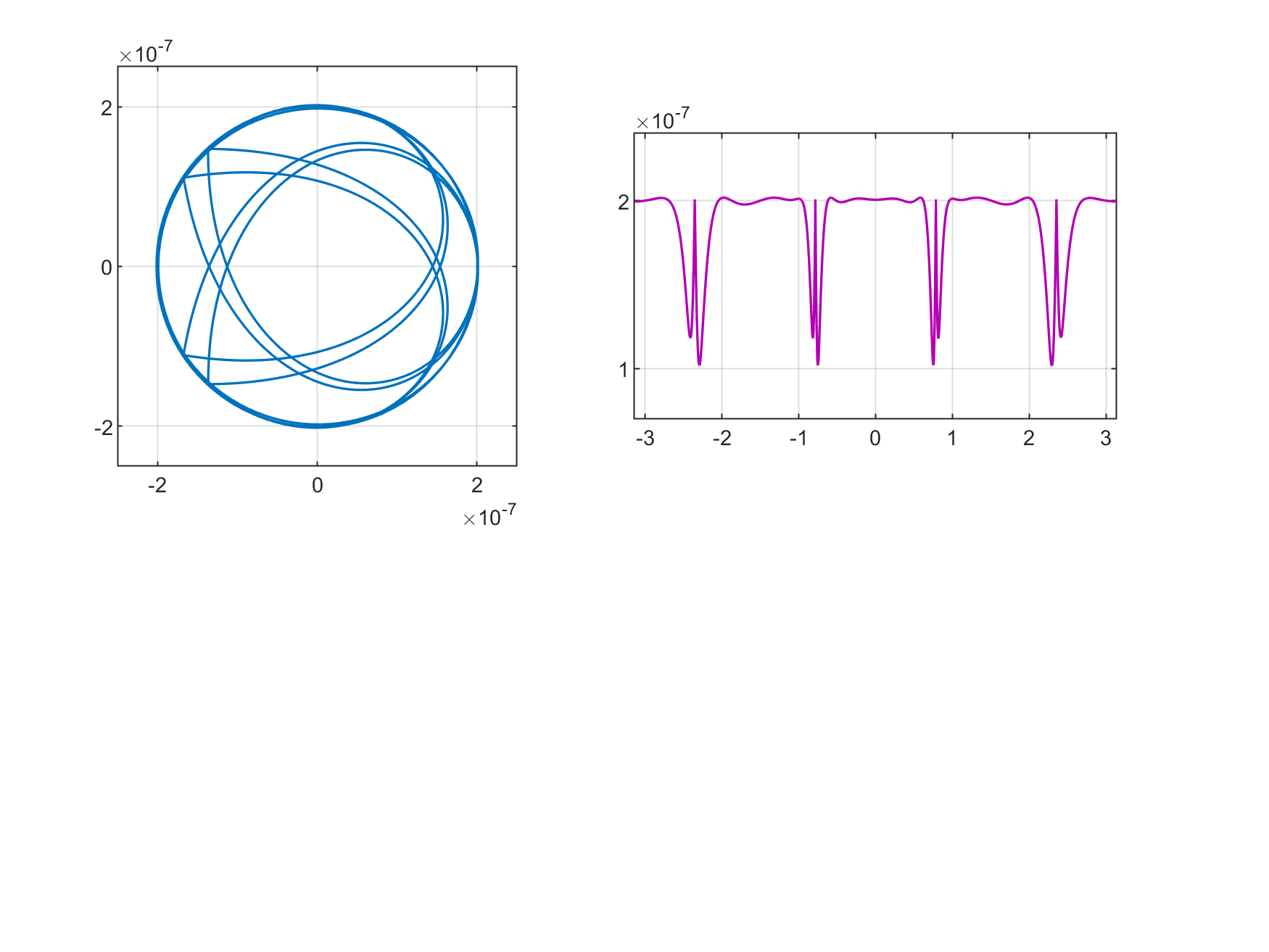}
\end{center}
\vskip -13pt
\caption{The image on the left shows the error curve $(f-r)(Z)$ for best degree $n=4$ approximation of $e^z$ on the unit
square by the AAA-Lawson algorithm~{\rm \cite{aaac,aaalawson}}.  The image on the
right shows $|(f-r)(z)|$ as a function
of $\hbox{\rm arg}(z)$.  The error curve is a curve
of winding number $2n+1 = 9$ with nearly constant modulus except near
the four corners~{\rm\cite{nearcirc}}.}
\end{figure}

\smallskip
{\em AAA approximation on a continuum.}  The original AAA algorithm works with a discrete
sample set $Z$, which in applications is typically the discretization of a continuum
by hundreds or thousands of points.  (Since poles and zeros will
cluster near singularities, it is important to 
cluster the sample points too near
corners of the domain.)  It is natural to want an algorithm that can deal with a
continuum more directly.   Recently a variant of AAA has
been developed for this and supported with open-source
MATLAB and Julia software~\cite{aaac}.  This adds speed, robustness,
and simplicity to a AAA computation.  For example, the command

{\small
\begin{verbatim}


         r = aaax(@abs,80);

\end{verbatim}
\par}

\noindent
computes a near-best degree 80 approximation
to $|x|$ on $[-1,1]$ in $0.2$~s on our laptop based on the
evaluation of $|x|$ in a total of only about 500 points.
The accuracy is about $3.8\times 10^{-10}$,
not too far from the best approximation error of this degree,
$4.4\times 10^{-12}$.  This is the classic problem made famous
by Donald Newman in 1964~\cite{newman}, a prototype of the general
phenomenon of root-exponential convergence of rational approximations
to functions with branch point singularities~\cite[chap.~25]{atap}.
It is a very difficult problem computationally, for the zeros and poles of $r$
cluster exponentially near the singularity at $x=0$ at distances as
small as $4.3\times 10^{-10}$~\cite{clustering}.
Varga, Carpenter, and Ruttan required many hours of computing in
200-digit extended precision arithmetic thirty years ago to compute
such approximations using the Remez algorithm with a $p/q$ representation~\cite{vrc}.
The currently most robust implementation of the Remez algorithm, the Chebfun
{\tt minimax} command, takes $100$~s to find this best approximation.

\smallskip
{\em AAA-LS algorithm for Laplace and related problems.}  Perhaps the
most fundamentally important and surprising extension and application of
AAA methods was introduced in an arXiv paper by Stefano Costa in 2020~\cite{costa1} and
is being actively extended in various directions~\cite{costa2,costa,xue}.  Since the
days of Joseph Walsh and Mstislav Keldysh nearly a century ago, it has
been recognized that if rational functions are good at approximating analytic
functions, then their real parts will be good at
approximating harmonic functions (i.e., solutions of the Laplace equation).
With AAA as a new tool for rational approximation, it was natural to
expect this to lead to new methods for solving Laplace and related problems.  However,
this is not straightforward.  For example, suppose we wish to solve 
$\Delta u = 0$ in a simply-connected
domain $\Omega$ subject to a real boundary condition $u(z)= h(z)$
on the boundary $\partial \Omega$.  AAA enables us quickly to calculate an
approximation $r\approx h$, but $r$ will have poles in $\Omega$ as well as outside.
In other words, it approximates an analytic extension of $h$ to a neighborhood of
$\partial\Omega$, but not to all of $\Omega$.  So how can it help in solving the
Laplace problem?

\begin{figure}
\begin{center}
\includegraphics[trim=60 10 60 0, clip, scale=.95]{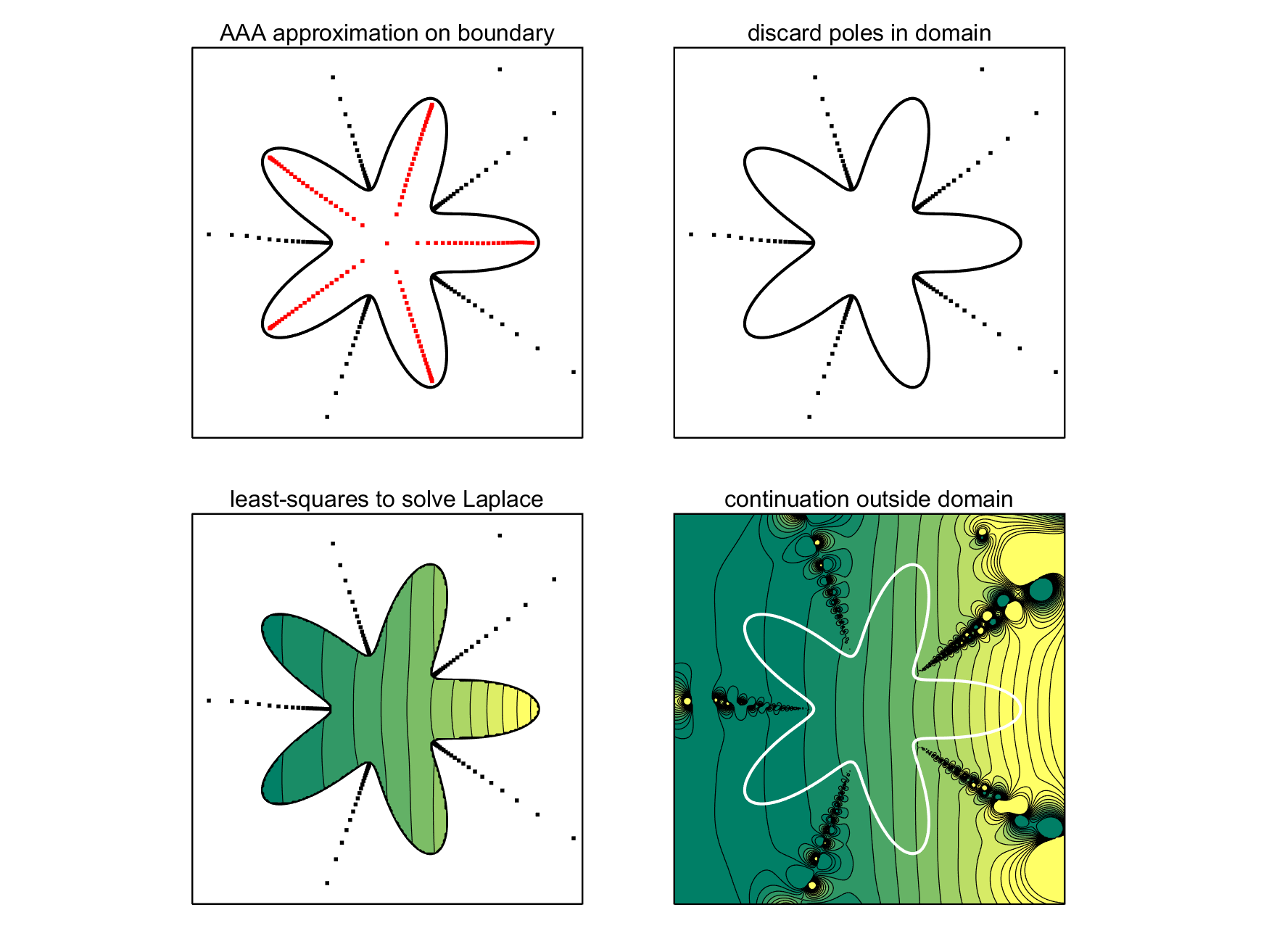}
\end{center}
\caption{The AAA-LS method applied to solve a Laplace Dirichlet
problem on a domain with five lobes. The solution computed in about $1$~{\rm s}, with
accuracy about\/ $7\times 10^{-7}$, is represented as the real part of
a rational function with $125$ poles and is evaluable in about
$1$ $\mu{\rm s}$ per point.  A polynomial representation of the same
accuracy would require a degree in the tens or hundreds of thousands.
The fourth panel shows evaluations outside
the problem domain, revealing branch cuts of the harmonic continuation
delineated by poles of the AAA approximation.}
\end{figure}

Costa realized that this challenge can be addressed by separating the
poles of the AAA approximant into those inside and those outside $\Omega$.
The functions $1/(z-s_j)$ defined by poles $\{s_j\}$ outside $\Omega$
span a space that contains excellent approximations to
the solution of the Laplace problem, which can be found numerically
by a linear least-squares calculation.  Thus the AAA-least squares
(AAA-LS) method consists of these steps, as illustrated
in Figure~3: (1) approximate the boundary
data by a rational function, (2) collect just the poles outside $\Omega$ to form
a basis for an approximation space,
(3) solve a least-squares problem in this basis.  Each pole leads to two
columns of the least-squares matrix, corresponding to the real and imaginary
parts of $1/(z-s_j)$.  The method proves extremely powerful in practice, as
is shown most recently in the computations of Y.~Xue of solutions to the biharmonic
equation for Stokes flow~\cite{xue} and 
in the solution of certain Poisson problems.
Initial theoretical justifications appear in~\cite{costa}, but there is much
more to be said, and further work on the theoretical side is underway. 

The separation of poles across two sides of a contour has many potential
applications beyond the solution of Laplace problems. This
is the starting point of effective numerical methods for Wiener--Hopf and
more general Riemann--Hilbert problems.  Closely related topics are the
Hilbert transform and Dirichlet-to-Neumann maps, and the AAA-LS approach 
is effective for all of these.  Examples can be found in~\cite{costa1,costa2,costa},
and there are many further possibilities.
In all these applications, rational functions offer the prospect of strikingly
fast and accurate computations in the presence of singularities or near-singularities,
when other methods tend to converge very slowly or require case-by-case analysis of
detailed behavior at singularities.

\smallskip
{\em Discussion.}  Unlike other algorithms for rational
approximation, such as the Ellacott--Williams and
Istace--Thiran algorithms or the well known method of Vector
Fitting~\cite{ew,vf,it}, AAA is not based on an attempt to enforce a
condition of optimality.  (Discussions of these and other algorithms
can be found in~\cite{aaa} and~\cite{aaalawson}.)  This appears to be one
of the secrets of its success.  As often happens in optimization,
the most difficult part of the problem may be not what to do near
the solution, but how to get near it in the first place.  AAA seems
extraordinarily good at this, for reasons not yet fully understood.
There is an interesting analogy with training algorithms for deep
learning, where optimization is carried out by methods such as
stochastic gradient descent that have little to do with optimality
conditions~\cite{goodfellow}.

As a generalization of polynomials to functions with poles
not constrained to lie at infinity, rational functions have a fundamental
importance in computational science, and AAA is opening many
doors.  Rational approximation, however, is only one long-established
nonlinear optimization problem among many.
Another important one is the approximation of a function $f(x)$ by a
linear combination of translates of Gaussians of unknown positions
and widths.  Before AAA, one would have imagined that these were
problems of comparable (and considerable) difficulty.
With the appearance of AAA, rational approximation suddenly seems
rather easy in practice.  Are there analogous algorithms for other
problems such as fitting by Gaussians?  We simply do not know,
and we regard this as a fascinating and important challenge.

\end{document}